\documentclass{article}

\usepackage{amsmath}
\usepackage{amssymb}
\usepackage{amsthm}
\usepackage{amsfonts}
\usepackage{amscd}
\usepackage{mathrsfs}
\theoremstyle{plain}
\newtheorem*{conjecture}{Conjecture}
\newtheorem*{main}{Main Theorem}
\newtheorem{theorem}{Theorem}[section]
\newtheorem{proposition}[theorem]{Proposition}
\newtheorem{corollary}[theorem]{Corollary}

\theoremstyle{definition}

\theoremstyle{remark}

\newtheorem{remark}[theorem]{Remark}
\newtheorem*{question}{Question}
\newtheorem{observation}[theorem]{Observation}

\numberwithin{equation}{section}

\def\R{ \mathbb R} 
\def\C{ \mathbb C} 
\def\Z{ \mathbb Z} 
\def\T{ \mathbb T} 
\def\P{ \mathbb P} 
\def\zmod#1{ {\mathbb Z}
/#1{\mathbb Z}}
\def\Hilb{\mathcal H}  
 
\def\d{ \partial }                  
\def\rttau{ \tau^{1/2} }        
\def\gotimes{\hat{\otimes}} 

\def\Tr{\text{\rm Tr}}
\def\Det{\text{\rm Det}}
\def\Pfaff{\text{\rm Pfaff}}

\def\Lie#1{ {\mathfrak #1} }

\def\ad{ \text{\rm ad}}
\def\fpairing#1#2{\langle #1 \wedge #2 \rangle } 
\def\pairing#1#2{\langle #1 , #2 \rangle }

\def\Line#1#2{ {\mathcal L}^{#1}#2 }   
\def\Mod#1#2{ {\mathcal M}_{#1}(#2) } 
\def\Cat#1#2{ {\mathcal C}_{#1}(#2) }  
\def\Bor#1{ {\mathscr B}{#1}}  
\def\Vec{ {\mathscr V}}  
\def\Mor#1#2{ {\mathcal G}_{#1}(#2) }  

\begin{document}
\title{Spin Chern-Simons and Spin TQFTs}
\author{Jerome A. Jenquin}
\date{May 9, 2006}
\maketitle
\begin{abstract}
In \cite{Je1} we constructed classical spin Chern-Simons for any compact Lie
group $G$: a gauge 
theory whose action depends on the spin structure of the 3-manifold.
Here we apply geometric quantization to the classical
Hamiltonian theory and investigate the formal properties of the partition
function in the Lagrangian theory, all in the case $G = SO_3$.  We find
that the quantum theory for $SO_3$ spin Chern-Simons corresponds to  the
spin TQFT constructed by Blanchet and Masbaum \cite{BM} in the same way 
that the quantum theory for standard $SU_2$ Chern-Simons corresponds to 
the TQFT constructed by Reshetikhen and Turaev \cite{RT} or the TQFT 
constructed by Blanchet, Habegger, Masbaum, and Vogel \cite{BHMV}.
\end{abstract}

 \section{Knot invariants and physics}\label{Introduction}
 
 The correlation between between knot invariants 
 and quantum Chern-Simons was first laid out in the 
 seminal paper by Witten \cite{Wi}.  The crux of Witten's
 argument is that the quantum theory is a 
 topological quantum field theory, or TQFT.  In 
 particular, he applies the formal properties of 
 the quantum partition function to show that the
 theory satisfies the defining axioms of a TQFT.
 Together with certain results from conformal
 field theory, Witten uses these axioms to 
 show that the quantum partition function 
 associated to a compact oriented 3-manifold
 is the Jones polynomial of the corresponding
 knot.  For a complete definition of a TQFT and
 an explanation of Witten's results we suggest \cite{A}.
 
 Mathematicians -- in particular, Reshetikhin and Turaev \cite{RT}--
 took a different approach to knot and 3-manifold
 invariants by constructing their own TQFT from the 
 representation theory of quantum groups.  So that 
 their invariants matched Witten's, their TQFT was 
 necessarily isomorphic to quantum Chern-Simons.
 
The quartet of Blanchet, Habegger, Masbaum, and Vogel (BHMV)
come at TQFTs from the other direction; their starting point is a 
generalized version of the Reshetikhin and Turaev invariant.  They apply
algebro-categorical techniques to the Kauffman 
bracket and an oriented bordism category \cite{BHMV}.  As 
desired, the BHMV TQFT nicely matches Witten's.  We elaborate in 
section \ref{BHMV review}.

Our interest here, though, is in the work by Blanchet and Masbaum that
followed the BHMV results.  In the same spirit, their starting
point is a refinement of the Reshetikhin and Turaev invariant --
one that now depends on the spin structure of the 3-manifold.  They 
apply the same algebro-categorical techniques, but now to a refined
version of the Kauffman bracket and a spin bordism category.
In this bordism category the manifolds and their boundaries must 
have compatible spin structures.  The resulting construction satisfies 
the properties of what Blanchet and Masbaum dub a ``spin TQFT",  
where the axioms that define a spin TQFT are essentially refined 
versions of the axioms that define the ``unspun" TQFTs mentioned above.

As we mentioned above, the BHMV TQFT corresponds well with
quantum Chern-Simons.  Upon learning of the Blanchet-Masbaum 
(BM) spin TQFT, one might naturally ask the following:

\begin{question} Is there a quantum field theory that
corresponds to the BM spin TQFT just as quantum Chern-Simons
corresponds to the BHMV TQFT; and if so, what is it?
\end{question}  

The answer , we claim, is yes; and the corresponding quantum field theory is the 
topic of this paper.  We refer to this field theory as spin Chern-Simons.

In the first section of this paper we review the 
relevant aspects of classical spin Chern-Simons .  These were 
worked out by the author in \cite{Je1} and we refer to that paper
for details and proofs.  In the second part we review the BHMV TQFT
and the BM spin TQFT.  Here we carefully point out the correspondence
between the BHMV TQFT and quantum Chern-Simons so that we can see how 
quantum spin Chern-Simons might correspond to the BM spin TQFT.  
We then address this hoped for correspondence on two fronts.  
The first is not necessarily rigorous.  We use the 
formal properties of the quantum partition function to show that
the 3-manifold invariants for both TQFTs display the same 
behavior.  On the second front we compute the dimensions
of the Hilbert spaces associated to closed, compact 2-manifolds
in quantum spin Chern-Simons.  We see that our dimension 
formulas are the same as those of the BM spin TQFT and that
the Hilbert spaces of spin Chern-Simons are refinements 
of the Hilbert spaces of standard Chern-Simons.  This in lock
step with the way in which the Hilbert spaces of the BM spin TQFT 
refine the Hilbert spaces of the BHMV TQFT.

\section{Classical spin Chern-Simons}

To define the classical theory we choose a compact 
Lie group $G$ and an orthogonal, rank zero virtual representation 
$\rho \in RO(G)$.  In practice, we often have an actual orthogonal
representation $\rho_0: G \rightarrow O(V)$ and take 
$$
\rho = \rho_0 - \dim V.
$$
That is, $\rho$ is often the difference of an actual representation
and the trivial representation of the same dimension.

We consider both the Lagrangian field theory over compact,
spun 3-manifolds, and the Hamiltonian theory over closed, 
compact, spun 2-manifolds.  In either case, the fields 
over a manifold $X$ are principal $G$-bundles with connection 
over $X$.  Technically, this ``space'' of fields is a category
$\Cat{G}{X}$ -- a groupoid, in fact.  The morphisms $\Mor{G}{X}$ for 
this category are the $G$-bundle isomorphisms that cover the identity on $X$; 
in other words, the gauge transformations.  Two
objects in $\Cat{G}{X}$ are considered physically equivalent 
if there is a gauge transformation taking one to the other.

\subsection{Classical Lagrangian theory}\label{Classical lagrangian}

Here, $X$ is a closed, compact, spun 3-manifold. 
The main object of interest in the Lagrangian field theory
is the action.  In our case, this is a $\T$-valued function 
on the space of fields $\Cat{G}{X}$, where $\T \subset \C$ 
are the unit modulus complex numbers .  To define the action
we place a Riemannian structure on $X$ and so obtain
a Dirac operator acting on sections of the spinor bundle
$$
D_X : \Gamma (S_X) \rightarrow \Gamma (S_X).
$$
From the pair $(P,A) \in \Cat{G}{X}$ -- where $P$ is a principal $G$-bundle 
over $X$ and $A$ is a connection on $P$ -- we obtain
the associated virtual vector bundle with connection $(\rho P, \rho A)$.  
We can couple $\rho A$ to the Dirac operator to obtain a twisted Dirac
operator
$$
D_X \otimes \rho A : 
\Gamma (S_X \otimes \rho P) \rightarrow 
\Gamma (S_X \otimes \rho P).
$$
This operator is elliptic, self-adjoint and quaternionic 
so that it has a discrete, real spectrum of eigenvalues with
even valued degeneracies.  This allows us to define our action
to be the assignment
\begin{align}
\Cat{G}{X} & \rightarrow \T  \notag\\
A & \mapsto \rttau (D_X \otimes \rho A)
\end{align}
where $\rttau (D)$ is a spectral invariant defined for 
any elliptic, self-adjoint, quaternionic operator $D$.
In fact, it is (a square root of) the exponentiated boundary term 
in the Atiyah-Patodi-Singer index theorem \cite{APS1} 
and we use this fact to our advantage in determining some 
of the action's properties.  In particular, if all of the data (the
bundle, connection, metric, and spin structure) bounds a 
principal G-bundle $P_M$ with connection $A_M$ over a spun 
4-manifold $M$ we can write

\begin{equation}\label{rttau: bounding case}
\rttau (D_X \otimes \rho A) = 
\exp \pi i \int_M \fpairing{\Omega^{A_M}}{\Omega^{A_M}}_{\rho} 
\end{equation}\noindent

where $\Omega^{A_M}$ is the curvature of the connection $A_M$ 
and

\begin{equation}\label{p1 form}
\pairing{\eta_1}{\eta_2}_{\rho} =
{-1 \over 8 \pi^2}  \Tr  \left( \rho (\eta_1) \rho (\eta_2) \right)
\end{equation}

is a bilinear Ad-invariant form on the Lie algebra $\Lie{g} = 
\text{Lie}(G)$.  Those familiar with standard Chern-Simons 
will note that the 4-dimensional integral is equal the action 
of that theory at the ``level'' determined by ${1 \over 2}\pairing{}{}$. 

We summarize other properties in the following theorems
whose proofs can be found throughout section (1) of \cite{Je1}
First we state how the action depends on the smooth 
parameters.
 
\begin{theorem}
Assume that $\dim \rho = 0$ and that $\pairing{}{}_{\rho}$
is a non-degenerate pairing.  Then $\rttau (D_X \otimes \rho A)$
is independent of the metric and invariant under gauge transformations
With respect to the connection $A$, the critical points are exactly the flat 
connections.  In other words, $d \rttau (D_X \otimes \rho A) = 0$ if and 
only if $\Omega^A =0$.
\end{theorem}   

Though we must choose a metric to define the action, the
choice is ultimately irrelevant and the physical theory does not
depend on it.  Also, we see that the classical fields, i.e. the critical points
of the action, are flat connections and that the action descends to the 
quotient $\Cat{G}{X} / \Mor{G}{X}$.  In these three respects spin
Chern-Simons is the same as standard Chern-Simons which is implicitly 
independent of the metric and whose gauge invariant action also has 
flat connections as its critical points.  

The next two theorems, in which we spell out the action's dependence on the 
discrete parameters, each require some preemptive explanation.  

One (of the two) discrete parameters is the spin structure
$\sigma$ on $X$.  As is well known, the equivalence classes
of spin structures is affine over $H^1(X; \zmod2)$.  Thus to
track how the action changes with respect to the shift $\sigma \mapsto
\sigma + \ell$ for some $\ell \in H^1(X; \zmod2)$ we consider the ratio
$$
q_{\sigma}(\rho P, \ell) =
{ \rttau (D_{\sigma + \ell} \otimes \rho A) \over
\rttau (D_{\sigma} \otimes \rho A)} .
$$  

\begin{theorem}\label{properties of q}
If $\sigma$, $\ell$, and $\rho$ are as above and $w_1(\rho P) = 0$ 
then 
$$
q_{\sigma}(\rho P, \ell) =
(-1)^{w_2(\rho P) \smile \ell}
$$
\end{theorem}
Thus a change in spin structure could cause the action to change
by a sign.  

The other parameter is the virtual representation $\rho \in RSO(G)$.
In standard Chern-Simons the analogue to $\rho$ is an Ad-invariant
bilinear form $\pairing{}{}$ on $\Lie{g}$.  Furthermore, in that theory, it is  
well known that the action's dependence on $\pairing{}{}$ factors through the 
Chern-Weil map
$$
\text{Sym}^2 \Lie{g}^* \rightarrow H^4(BG; \R).
$$
Similarly, in spin Chern-Simons, we have the following.

\begin{theorem}
The action's dependence on
$\rho$ factors through a homomorphism
$$
\lambda_G : RSO(G) \rightarrow E^4(BG)
$$
where $E^{\bullet}(\cdot)$ is a generalized cohomology. 
\end{theorem}
The generalized cohomology in question has a
corresponding spectrum of spaces $\{ E^j \}$ such that 
$E^j(\cdot) = [\cdot , E^j]$.  In particular, the space $E^4$ is 
(up to homotopy) the bottom two rungs of the Postnikov
tower for $BSO$.  We will not say more here about $E^4(BG)$
in general but defer to the appendix in \cite{Je1}. 

However, because this paper focuses on the cases $G = SU_2$
and $G = SO_3$, we do consider the following examples.
Let $\rho_0 : SU_2 \rightarrow SO_4$ denote the realization of the
standard action on $\C^2 = \R^4$ and let $id_{SO_3}: SO_3 \rightarrow
SO_3$ denote the standard action on $\R^3$.  With some foresight 
we define 
$$
{\bf 1'} = \lambda_{SU_2} (\rho_0 - 4) \quad \quad
\text{and} \quad \quad
{\bf 1} = \lambda_{SO_3} (\text{id}_{SO_3} - 3).
$$  
Then it turns out that
$$
E^4(BSU_2) = \Z \cdot {\bf 1'} \quad \quad
\text{and} \quad \quad
E^4(BSO_3) = \Z \cdot {\bf 1}.
$$
and that the covering homomorphism $SU_2 \rightarrow SO_3$ 
induces the homomorpshim

\begin{align}
E^4(BSO_3) & \rightarrow E^4(BSU_2) \notag \\
n \cdot {\bf 1} & \mapsto 2n \cdot {\bf 1'}.
\end{align} 

This fact plays a crucial role in determining the formal
properties of the quantum partition function and in relating
the Hamiltonian theories for $SU_2$ and $SO_3$.

\subsection{Classical Hamiltonian theory}\label{Classical hamiltonian}

Here $Y$ is a closed, compact, spun 2-manifold. 
The main object of interest in the Hamiltonian field theory
is the prequantum line bundle over the classical phase space.
To define it, we place a Riemannian structure on $Y$ and so obtain
a chiral Dirac operator 
$$
D_Y : \Gamma (S^+_Y) \rightarrow \Gamma (S^-_Y).
$$
Much as before, to any pair $(P,A) \in \Cat{G}{Y}$  
we associate a twisted chiral Dirac operator
$$
D_Y \otimes \rho A : 
\Gamma (S^+_Y \otimes \rho P) \rightarrow 
\Gamma (S^-_Y \otimes \rho P).
$$
This operator is elliptic and skew-symmetric so that we
can define the assignment
 \begin{align}\label{Pfaff}
A & \mapsto \Pfaff^{-1}(D_Y \otimes \rho A ) = 
\bigwedge\nolimits^{\text{top}} \ker (D_Y \otimes \rho A)
\end{align}
where $\Pfaff^{-1}(D)$ is the {\it inverse Pfaffian} line associated 
to a skew-symmetric operator $D$.  If we fix a $G$-bundle
$P \rightarrow Y$, then these lines fit together to form 
a smooth line bundle 
$$
\Line{\rho}{(P)} \rightarrow \Cat{}{P} 
$$ 
where $\Cat{}{P}$ denotes the (affine) space of connections on $P$.
In fact, this line bundle has a natural hermitian structure
and compatible connection \cite{F2}.

We mention two important aspects  of this inverse
Pfaffian line bundle.  The first is in regards to the holonomy
associated to the natural connection.   If 
$\gamma: S^1 \rightarrow \Cat{}{P}$ is a closed path of $G$-connections
over $Y$, it induces a $G$-connection $A_{\gamma}$ over the 
3-manifold $S^1 \times Y$.  Then the holonomy around $\gamma$
is given by \cite{BF2} 
\footnote{Generally, the holonomy is given by
taking an adiabatic limit over metrics on $S^1$.  In this case, metric
independence eliminates the need to take this limit.}

\begin{equation}\label{Holonomy}
\text{hol}_{\gamma} = \rttau (D_{S^1 \times Y} \otimes \rho A_{\gamma}).
\end{equation}
   
The second aspect is in regards to how $G$-connection automorphisms
lift to the bundle.  If $\phi$ is a gauge transformation that preserves the
$G$-connection $A$, then $\phi$ induces an automorphism of 
the line $\Pfaff^{-1}(D_Y \otimes \rho A)$.  To compute that
automorphism, we note that $\phi$  and $A$ induce 
another $G$-connection $A_{\phi}$ over $S^1 \times Y$.  Then the
induced automorphism is given by \cite{Je1}

\begin{equation}\label{Automorphism}
\text{aut}_{\phi} = \rttau (D_{S^1 \times Y} \otimes \rho A_{\phi}).
\end{equation}

In both cases the spin structure on $S^1 \times Y$ is the product spin 
structure induced by the given spin structure on $Y$ and the {\it bounding}
spin structure on $S^1$.  This is the spin structure that extends to the
disc.

We now address the fact that we had to choose a metric on $Y$ to
define the inverse Pfaffian bundle.
The line bundle is independent of the metric
on $Y$ in the sense that, given a different metric, the two 
line bundles are canonically isomorphic.  This follows immediately
from the metric independence of the action and \ref{Holonomy}.
Thus, in this sense, the Hamiltonian theory for spin Chern-Simons is 
independent of the metric much like it implicitly is for standard Chern-Simons.  

Just as it is on compact 3-manifolds, the space of classical 
solutions is still the category of flat $G$-connections over $Y$.
Also, we still consider two $G$-connections to be equivalent if they 
lie in the same $\Mor{G}{Y}$-orbit.  Given \ref{Pfaff} it is clear that any gauge 
transformation between two objects of $\Cat{G}{Y}$ induces a natural 
isomorphism between their corresponding inverse Pfaffian lines.  Altogether 
this gives us the 
{\it prequantum line bundle} 
$$
\Line{\rho}{(Y)} \rightarrow \Mod{G}{Y}
$$
where $\Mod{G}{Y}$ is the moduli stack of flat $G$-bundles 
over $Y$.  The claim is that, over the flat $G$-connections, the 
inverse Pfaffian bundle described above, along with all of its
geometry, descends to the quotient under gauge transformations.
This is proven in \cite{Je1}.  Here we have denoted the descendant 
line bundle by $\Line{\rho}{(Y)}$.  As mentioned above, the prequanutm
line bundle over $\Mod{G}{Y}$ is the main object of interest in the 
classical Hamiltonian theory.  It also plays a large role in the quantum
theory, as we see in section \ref{quantum hamiltonian}.  

Before we end this review of classical spin Chern-Simons we point out
one final property of the prequantum line bundle.  The fibers of the bundle
$\Line{\rho}{(Y)}$ each have a natural $\zmod2$ grading given by 
$$\label{Grading}
| \Line{\rho}{(Y)}_{A} | = \dim \ker (D_Y \otimes \rho A) \pmod2 .
$$  
In other words, the grading is given by the mod-2 index of the 
skew-adjoint operator $D_Y \otimes \rho A$.  As this is a topological
invariant the grading is locally constant over the moduli stack.  This 
grading plays an important role in the quantum theory.

This ends our review of classical spin Chern-Simons.

\section{The BHMV TQFT and BM spin TQFT}

To better understand the correspondence between quantum spin Chern-Simons
and the BM spin TQFT we first review the correspondence between the BHMV 
TQFTs and quantum Chern-Simons.  We also review the relevant aspects of the 
BM spin TQFT.  See \cite{BHMV} and \cite{BM} respectively for details of the 
constructions and more details regarding the results.

\subsection{The BHMV TQFT versus the Witten TQFT}\label{BHMV review}

The BHMV TQFTs are constructed using combinatorial-topological techniques
and categorical machinery in conjunction with the knot-theoretic {\it Kauffman 
bracket}.  The result is a family of functors $V_p : \Bor{} \rightarrow \Vec$ 
indexed by positive integers $p \in \Z^{>0}$.  The objects of the domain
category $\Bor{}$ are closed, oriented 2-manifolds and its morphisms
are 3-manifolds with boundary.  
\footnote{The domain category considered in \cite{BHMV} is really that of 
2-manifolds with {\it $p_1$ structure}, which is a central extension of the 
category $\Bor{}$.  For simplicity we demur the issue of $p_1$ structures.}  
Thus if $Y_j$, $j=1,2$, are two objects of $\Bor{}$ and $X$ is a 3-manifold 
such that $\d X = -Y_1 \sqcup Y_2$, then $X$ is a morphism $Y_1 
\rightarrow Y_2$.  The objects of the codomain category $\Vec$ are finite 
dimensional, complex vector spaces and its morphisms are complex linear maps.  
\footnote{The codomain category really considered in \cite{BHMV} are 
modules over an abstract cyclotomic field $k_p$.  We have taken the liberty of 
choosing a particular extension of $k_p$ to $\C$, one for which the 3-manifold 
invariants of the Witten and BHMV theories agree.}
Thus, in the example above, $V_p(X)$ is a morphism $V_p(Y_1) \rightarrow
V_p(Y_2)$.

To the empty 2-manifold the functors assign $V_p(\emptyset) = \C$.  If,
for example, $\d X = Y$ then $X$ is a morphism $\emptyset \rightarrow Y$,
so that $V_p(X)$ is a morphism $\C \rightarrow V_p(Y)$; or what is the
same, $V_p(X)$ is a vector in $V_p(Y)$.  If, instead, $\d X = \emptyset$ then 
$X$ is a morphism $\emptyset \rightarrow \emptyset$ so that $V_p(X)$ is
a morphism $\C \rightarrow \C$; or what is the same, $V_p(X)$ is an element
of $\C$.  For any object $Y$ of $\Bor{}$, the vector space $V_p(Y)$ has a 
non-degenerate hermitian inner product $\pairing{}{}_Y$ such that if $\d X_1 
= \d X_2 = Y$, then
$$
\pairing{V_p(X_1)}{V_p(X_2)}_Y =
V_p (X_1 \cup_Y (-X_2)) \in \C.
$$ 
On top of that, there are natural isomorphisms $V_p(-Y) \rightarrow 
\overline{V_p(Y)}$ and $V_p(Y_1) \otimes V_(Y_2) \rightarrow
V_p(Y_1 \bigsqcup Y_2)$.  Thus the functors $V_p$ satisfy the axioms
of a 2-dimensional topological quantum field theory \cite{A}.

As is well known, physicists believe that classical $SU_2$ Chern-Simons, 
in conjunction with the Feynmann path integral, defines a similar family 
of functors $Z_k : \Bor{} \rightarrow \Vec$ \cite{Wi}.  The indexing 
set consists of positive elements $k \in H^4(BSU_2) \cong \Z$.  The formal
properties of the Feynmann path integral offer an easy ``proof'' that the 
functors $Z_k$ each define a 2-dimensional TQFT.  This is assuming,
of course, the one accepts the path integral.  Taking this for 
granted, we can make the following observation.   

\begin{observation}\label{obs1}
If one compares the functors $Z_k$ and $V_p$ then, whenever $p = 2(k+2)$,
one sees that 
\begin{itemize}
\item
For any object $Y$ of $\Bor{}$ the hermitian vector
spaces $V_{2(k+2)}(Y)$ and $Z_k(Y)$ have the same dimension.
In particular, if $Y$ is a genus $g$ 2-manifold then
\begin{equation}\label{verl form}
\dim V_{2(k + 2)}(Y) = \dim Z_k(Y) =
\left( \frac{k+2}{2} \right)^{g-1} 
\sum_{j=1}^{k+1} 
\left( \sin \frac{\pi j}{k+2} \right)^{2-2g} ,
\end{equation} 
which is the famous Verlinde formula \cite{Ve}.
\item
For many examples of closed, oriented 3-manifolds $X$, 
it has been shown that $V_{2(k+2)}(X) = Z_{k}(X)$ 
(cf. \cite{FG} , \cite{Wi}, and \cite{Li1}).
\item   
So that we can point out one last feature, we must first recall a 
fact about the algebraic topology of closed compact 2-manifolds.
For any object $Y$ of $\Bor{}$ there is a central extension
of $H^1(Y; \zmod2)$ by $\zmod4$; or what is
the same, there exists a long exact sequence
$$
0 \rightarrow \zmod4 \longrightarrow \Gamma (Y) 
\longrightarrow H^1(Y; \zmod2) \rightarrow 0.
$$
$\Gamma (Y)$ is a quotient of the usual Heisenberg group
$H (Y)$ which itself fits into the short exact sequence 
$$
0 \rightarrow \Z \longrightarrow H (Y) 
\longrightarrow H^1(Y; \Z) \rightarrow 0.
$$
The last feature we wish to point out it that the vector spaces
$V_{2(k+2)}(Y)$ and $Z_k(Y)$ each support a natural 
$\Gamma (Y)$ action and are equivalent as representations
of $\Gamma (Y)$ \cite{AM}.
\end{itemize}
\end{observation}

We make some remarks about this last observation as it has
bearing on the $SO_3$ spin-Chern-Simons theory.  At the positive
$SU_2$ levels $k$ the vector spaces 
$V_{2(k+2)}(Y)$ and $Z_k(Y)$ decompose as representations of 
$\Gamma (Y)$.  In the cases $k \equiv 0$ (mod 4) there is a natural one-one 
correspondence between irreducible components and  
$\zmod2$-bundles on $Y$; and in 
the cases $k \equiv 2$ (mod 4) there is a natural one-one correspondence 
between irreducible components and spin structures on $Y$.  The 
latter correspondence will appear again when we consider the spin 
TQFTs of Blanchet and Masbaum in the next subsection and then yet 
again in the quantum spin-Chern-Simons theory for $SO_3$.  In fact, 
we believe that the former and latter correspondences are part of a 
larger $SO_3$ gauge theory whose consideration we reserve for another 
time.  In the cases $k \equiv 1$ or $3$ (mod 4) the vector spaces 
decompose but the irreducible components do not correspond to any 
topological structures.   

As a clue to why the 0 (mod 2)-valued $SU_2$ levels should
have anything to do with an $SO_3$ spin-Chern-Simons theory,
we recall that levels of the later are elements of $E^4(BSO_3)$.
At the end of section \ref{Classical lagrangian} we stated that $E^4(BSO_3) =
 \Z  \cdot {\bf 1}$ and that $E^4(BSU_2) = \Z \cdot {\bf 1'}$.  Furthmore we
stated that the $2:1$ covering map $\beta : SU_2 \rightarrow SO_3$
induces a homomorphism 
\begin{align}
\beta^* : E^4(BSO_3) &\longrightarrow E^4(BSU_2) \\ \nonumber
              {\bf 1} & \longmapsto  2 \cdot {\bf 1'}
\end{align}            
Thus, the 0 (mod 2)-valued $SU_2$ levels correspond to
$SO_3$ levels.  We will have much to say about this is the sections
to follow.  For now we collate all of these correspondences in the 
following table: 

\vskip 1 cm
\begin{center}
Corresponding Representations of $\Gamma (Y)$
\end{center}

\begin{tabular}[c]{| c | c | c | c | c | }
\hline
BHMV level (mod 8) & 0 & 4  & 2  & 6   \\
\hline
$SU_2$ level (mod 4)  & 2  & 0 & 1 & 3  \\
\hline
$SO_3$ level (mod 2)  & 1 & 0 &   &  \\
\hline
Topological Structure  & spin structure & $\zmod2$-bundle & &  \\
\hline  
\end{tabular}

\vskip 1 cm
Before moving on, we mention that our considerations throughout the
rest of this paper focus on the corresponding representations 
in first column; i.e. those representations that correspond to spin structures
on the 2-manifold. 

\subsection{The BM spin TQFT}\label{BM review}

The BM spin TQFTs are constructed using the same combinatorial-topological 
techniques and categorical machinery as was used in constructing the BHMV 
TQFTs; but they are used in conjunction with a knot-theoretic invariant that
is sensitive to spin structures.  The result is a family of functors $V^s_p : 
\Bor{}^s \rightarrow \Vec^s$ indexed by positive 0 (mod 8)-valued integers 
$p \in 8\Z^{>0}$.  The objects of the domain category $\Bor{}^s$ are 
closed, spin 2-manifolds and its morphisms are spin 3-manifolds with boundary.  
\footnote{Again, the domain category considered in \cite{BM} is really that 
of spin 2-manifolds with $p_1$ structure; or what is the same, 2-manifolds
with {\it string stucture}.}  
Thus if $(Y_j, \sigma_j)$, $j=1,2$, are two objects of $\Bor{}^s$ and 
$(X, \Sigma)$ is a spin 3-manifold such that $\d (X,\Sigma) = (-Y_1 \sqcup 
Y_2, -\sigma_1 \sqcup \sigma_2)$, then $(X,\Sigma)$ is a morphism 
$(Y_1,\sigma_1) \rightarrow (Y_2,\sigma_2)$.  The objects of the codomain 
category $\Vec^s$ are finite dimensional, $\zmod2$-graded complex vector spaces 
and its morphisms are complex linear maps that preserve the grading.
\footnote{Again, the codomain category really considered in \cite{BM} is
that of $\zmod2$-graded modules over an abstract cyclotomic field $k_p$.  
We take the same liberties of we took for the BHMV TQFT in extending $k_p$ to 
$\C$.}

All the axioms satisfied by the BHMV TQFTs are also satisfied by the BM spin
TQFTs but with one caveat.  For any object $(Y, \sigma)$ of 
$\Bor{}^s$ we can write 
$$
V^s_p(Y, \sigma) = V^s_{p,0}(Y, \sigma) \oplus
V^s_{p,1}(Y, \sigma)
$$
where $V^s_{p,0}(Y,\sigma)$, $V^s_{p,1}(Y,\sigma)$ are respectively 
the even and odd components of $V^s_p(Y,\sigma)$.  Then the caveat is that
there are natural isomorphisms
\begin{align}\label{gr1}
& V^s_{p,0}(Y_1 \sqcup Y_2, \sigma_1 \sqcup \sigma_2)
\rightarrow \\
V^s_{p,0}(Y_1, \sigma_1) & \otimes V^s_{p,0}(Y_2, \sigma_2) \oplus
V^s_{p,1}(Y_1, \sigma_1) \otimes V^s_{p,1}(Y_2, \sigma_2) \nonumber
\end{align}
and
\begin{align}\label{gr2}
& V^s_{p,1}(Y_1 \sqcup Y_2, \sigma_1 \sqcup \sigma_2)
\rightarrow \\
V^s_{p,0}(Y_1, \sigma_1) & \otimes V^s_{p,1}(Y_2, \sigma_2) \oplus
V^s_{p,1}(Y_1, \sigma_1) \otimes V^s_{p,0}(Y_2, \sigma_2). \nonumber
\end{align}
This is expressed more succinctly in terms of the graded tensor product
``$\gotimes$''.  Indeed, in the language of $\zmod2$-graded vector 
spaces \eqref{gr1} and \eqref{gr2} amount to saying that there is a
natural isomorphism 
$$
V^s_{p}(Y_1 \sqcup Y_2, \sigma_1 \sqcup \sigma_2)
\rightarrow
V^s_{p}(Y_1, \sigma_1) \gotimes V^s_{p}(Y_2, \sigma_2). 
$$ 

We enumerate some of the more interesting features of the BM spin 
TQFTs: 

\begin{enumerate}
\item[1.]
For a closed 3-manifold $X$ the BM invariants $V^s_p(X,\sigma)$ are 
refinements of the BHMV invariant $V_p(X)$ in the sense that 
$\sum_{\sigma} V^s_p(X,\sigma) = V_p(X)$.
\item[2a.]
For a closed, genus $g$, spin 2-manifold $(Y,\sigma)$
\begin{equation}\label{BM verl even}
\dim V^s_{p,0}(Y, \sigma) =
\frac{1}{2^{2g}}(\dim V_p(Y) + 
(\frac{p}{4})^{g-1}((-1)^{\epsilon(\sigma)}2^g - 1))
\end{equation}
where $\epsilon (\sigma)$ is the {\it Arf invariant} of the 
spin structure.
\item[2b.]
For the same spin 2-manifold
\begin{equation}\label{BM verl odd}
\dim V^s_{p,1}(Y, \sigma) =
\frac{1}{2^{2g}}(\dim V'_p(Y) - 
(\frac{p}{4})^{g-1}((-1)^{\epsilon(\sigma)}2^g - 1))
\end{equation}
where $V'_p(Y)$ is a particular vector space associated to $Y$ by 
the BHMV TQFT (see Remark 5.11 of \cite{BHMV}).  
At any rate its dimension is given by
\begin{equation}\label{verl odd}
\dim V'_p(Y) = \left( \frac{p}{4} \right)^{g-1}
\sum_{j=1}^{p/2 -1} (-1)^{j+1} 
\left( \sin \frac{2 \pi j}{p} \right)^{2-2g}.
\end{equation}
\item[3.]
There exists a canonical isomorphism
$$
V_p(Y) \longrightarrow \bigoplus_{\sigma}
V^s_{p,0}(Y,\sigma).
$$
\end{enumerate}

\subsection{A conjecture and the Main Theorem}

The points of Observation \ref{obs1} amount to evidence toward
a conjecture; the conjecture being that the quantum $SU_2$ Chern-Simons
theory is, in fact, a TQFT and that the $SU_2$ TQFT at level $k$ is 
isomorphic to the BHMV TQFT at level $p = 2(k+2)$.  That this is 
{\it only} a conjecture comes from the fact that the Feynmann path integral 
is, at this point, not a well-defined mathematical object.  In fact, it is very
likely a conjecture that will not be proven any time soon.  Despite that,
the evidence offered in Obervation \ref{obs1} does make for a rather
convincing empirical argument.

If the conjecture is to be believed it seems that the knot-theoretic BHMV
TQFTs have a gauge-theoretic correspondence in the quantum $SU_2$ 
Chern-Simons theory.  Upon discovering the knot-theoretic BM spin TQFTs one is 
naturally lead to ask if they have their own gauge-theoretic correspondence.
This is the question that motivates our investigation of the quantum $SO_3$ 
spin-Chern-Simons theory.  Indeed, we offer the following conjecture.

\begin{conjecture}
The quantum $SO_3$ spin-Chern-Simons theory is, in fact, a TQFT and
at the 1 (mod 2)-valued $SO_3$ level $2m-1$, it is isomorphic to 
the BM spin TQFT at level $p = 8m$.
\end{conjecture}  

Physicists would like to believe that the quantum partition function
and canonical quantization can generate a spin TQFT from the classical
theory described above.  While it is (at this point) impossible to rigorously
do so we can, nonetheless, exploit the features of this quantization map to
compare the conjectured spin TQFT and the BM spin TQFT.  In particular,
we appeal to the formal properties of the quantum partition function 
and the procedure of {\it geometric quantization}. 

Now consider spin Chern-Simons, as defined above, for a 
compact group $G$ and level $k$.  Let $Z_{G}(X,\sigma, k )$ denote 
the quantum partition function associated to a closed spin 3-manifold 
$(X, \sigma)$, which we consider only formally.  And let $\Hilb_{G}(Y,\sigma, k)$ 
denote the Hilbert space associated to  a closed spin 2-manifold $(Y, \sigma)$
via geometric quantization.  As evidence toward the conjecture above we provide 
the following theorem.  

\begin{main}
For odd-valued $SO_3$-levels $k = 2m-1$ we have the following:
\begin{enumerate}
\item[1]
{\rm $Z_{SO_3,2m-1}(X^3,\sigma, )$ has the same 
formal properties as $V^s_{8m}(X^3,\sigma )$ in the sense that }
$$
Z_{SU_2, 4m-2} = \sum_{\sigma}
Z_{SO_3,2m-1}(X^3,\sigma, ).
$$ 
\item[2]
{\rm $\Hilb_{SO_3,2m-1}(Y^2,\sigma )$ is $ \Z /2\Z $ graded.  That is,}
$$
\Hilb_{SO_3,2m-1}(Y^2,\sigma ) = \Hilb_{SO_3,2m-1}^0(Y^2,\sigma ) 
\oplus \Hilb_{SO_3,2m-1}^1(Y^2,\sigma ).
$$
\item[3.a]
$\dim \Hilb_{SO_3,2m-1}^0(Y^2,\sigma ) =
\dim V^s_{8m,0}(Y^2,\sigma ) \quad$ {\rm and}
\item[3.b]
$\dim \Hilb_{SO_3,k}^1(Y^2,\sigma ) =
\dim V^s_{8m,1}(Y^2,\sigma ).$
\item[4]
$\Hilb_{SU_2, 4m-2}(Y^2) = \bigoplus_{\sigma}
\Hilb_{SO_3,2m-1}^0(Y^2,\sigma ).$
\end{enumerate}
\end{main}

\section{Formal Properties of the Three-Manifold Invariants}

\subsection{Spin Chern-Simons versus Chern-Simons}
In this short section we show that if we consider the quantum partition 
function as a formal object then we can show that the spin 3-manifold
invariants for the $SO_3$ spin-Chern-Simons theory at $SO_3$ level
$k = 2m-1$ are refinements of the 3-manifold invariants for the 
Chern-Simons theory at $SU_2$ level $k' = 4m-2$.

To begin we prove the following proposition.

\begin{proposition}\label{simply conn}
Let $\rho$ be a real rank zero virtual represenation of a connected, simply 
connected, compact Lie group $G$ and let $\pairing{}{}_{\rho}$ denote
the symmetric pairing defined by \eqref{p1 form}.  If $A$ is a 
$G$-connection over a closed, spin 3-manifold $X$, we let
$\exp 2 \pi i S_X (A)$ denote the $\T$-valued Chern-Simons 
invariant determined by the pairing $\frac{1}{2}\pairing{}{}_{\rho}$, 
as defined in \cite{F2}.  Then 
$$
\rttau_X (D_{\rho A}) = \exp 2 \pi i S_X (A).
$$ 
\end{proposition} 

\begin{proof}
The proof relies on the fact that the cobordism group $\Omega_3^{spin}(BG) = 0$
whenever $G$ is compact and simply connected.  In that case there exists a
spin 4-manifold $M$ such that $\d M = X$ as a spin manifold and there exists
an extension $A'$ of $A$ over $M$.  On the one hand it is well-know that
$$
\exp 2 \pi i S_X (\rho A) = 
\exp 2 \pi i \int_M \frac{1}{2}
\pairing{\Omega^{A'}}{\Omega^{A'}}_{\rho}.  
$$
On the other hand, since $\rho$ has rank zero the APS index theorem
implies
$$
\rttau_X (D_{\rho A}) =
\exp \pi i \int_M  
\pairing{\Omega^{A'}}{\Omega^{A'}}_{\rho},
$$
and this proves the proposition
\end{proof}

To make the equality of these theories even stronger we point out the
correspondence between their respective levels.
Recall that for Chern-Simons the levels are elements of $H^4(BG)$
while for spin-Chern-Simons the level are elements of $E^4(BG)$.
However, for $G$ simply connected there is a natural isomorphism
$i : H^4(BG) \rightarrow E^4(BG)$ so that Chern-Simons theory
at level $\alpha \in H^4(BG)$ is isomorphic to spin Chern-Simons
at level $i (\alpha) \in E^4(BG)$.  In particular, $SU_2$ Chern-Simons
at $SU_2$ level $k \in H^4(BSU_2) \cong \Z$ is isomorphic to 
$SU_2$ spin-Chern-Simons at level $k \in E^4(BSU_2) \cong \Z$.

\subsection{$SU_2$ and $SO_3$-connections}

Let $X$ be a closed spin 3-manifold.  The 2nd Steifel-Whitney class provides 
a one-one correspondence between isomorphism classes of $SO_3$ bundles
on $X$ and elements of $H^2(X;\zmod2)$.  Recall that $\overline{\Cat{G}{X}}$ 
denotes the quotient space of all $G$-connections on $X$ with respect to the gauge 
group $\Mor{G}{X}$.  Then the components of the quotient
space $\overline{\Cat{SO_3}{X}}$ are indexed by elements of $H^2(X;\zmod2)$
and we let $\Cat{SO_3, b}{X}$ denote the category of $SO_3$-connections
$(P,A)$ such that $w_2(P) = b$.

All $SU_2$ bundles over $X$ have a section so that the quotient space 
$\overline{\Cat{SU_2}{X}}$ has only one component.  The standard $2:1$
covering homomorphism induces a functor $\beta : 
\Cat{SU_2}{X} \rightarrow \Cat{SO_3, 0}{X}$.  A simple
argument shows that if two $SU_2$-connections map to the same 
$SO_3$-connection then they differ by a unique $\zmod2$ bundle
(up to a global choice of sign); and if two $SU_2$-connections differ 
by a $\zmod2$ bundle then they map to the same $SO_3$-connection.
Thus 
$$
\overline{\Cat{SO_3,0}{X}} =
\overline{\Cat{SU_2}{X}} / \overline{\Cat{\zmod2}{X}}.
$$  

From the discussion above it would seem that, at least formally,
\begin{align}\label{formal int1}
&\int_{\overline{\Cat{SU_2}{X}}} \quad \bullet \quad \mu_{SU_2}(X) = \\
\int_{\overline{\Cat{SO_3,0}{X}}} & 
\int_{\overline{\Cat{SU_2}{X}}/ \overline{\Cat{SO_3,0}{X}}} \quad 
\bullet \quad \mu_{S0_3,0}(X) \cdot \mu_{\zmod2}(X) \nonumber
\end{align}
where $\mu_G(X)$ is meant to be some sort of ``measure'' on $\overline{\Cat{G}{X}}$.
We point out that, while integrating over $\overline{\Cat{SU_2}{X}}$ or 
$\overline{\Cat{SO_3}{X}}$ does not make sense mathematically, integrating over
the finite fibers of $\overline{\Cat{SU_2}{X}} \rightarrow 
\overline{\Cat{SO_3,0}{X}}$ involves a finite sum that does make sense.

\subsection{Some formal manipulations}

As an example of a formal manipulation using \ref{formal int1},
imagine that $f : \overline{\Cat{SO_3}{X}} \rightarrow \C$ is a measurable
function.  Then, according to \ref{formal int1},
\begin{align}\label{formal int2}
\int_{\overline{\Cat{SU_2}{X}}} \quad f \circ 
\beta \quad \mu_{SU_2}(X) & = \\ \nonumber
\int_{\overline{\Cat{SO_3,0}{X}}} 
\int_{\overline{\Cat{SU_2}{X}}/ \overline{\Cat{SO_3,0}{X}}}\quad f \quad
\mu_{S0_3,0}(X) \cdot \mu_{\zmod2}(X) & = \\ 
\#(H^1(X;\zmod2)) \int_{\overline{\Cat{SO_3,0}{X}}} \quad f \quad 
\mu_{S0_3,0}(X) & \nonumber 
\end{align}
This is exactly the situation we encounter below.

On the one hand we consider $SO_3$ spin Chern-Simons at level $2m-1$.  
The action is given by 
$$
A \longmapsto f_{\sigma}(A) = 
\rttau_{X,\sigma} (D_{\rho A})^{2m-1}.  
$$    
where $\rho = \text{id} - 3$ is the standard representation minus the 
3-dimensional trivial representation and $\sigma$ is the spin structure on $X$.  
On the other hand we consider $SU_2$ Chern-Simons at level $4m - 2$ in 
which case a simple argument shows the action is given by 
\begin{equation}\label{def f'}
A' \longmapsto f'(A') = f_{\sigma} \circ \beta (A').
\end{equation}
We consider the formal integrals
\begin{align*}
&Z_{SO_3}(X,\sigma, 2m-1) = \int_{\overline{\Cat{SO_3}{X}}}
\quad f_{\sigma} \quad \mu_{S0_3}(X) 
\quad \text{and} \\
&Z_{SU_2}(X, 4m-2) = \int_{\overline{\Cat{SU_2}{X}}} \quad 
f' \quad \mu_{SU_2}(X), 
\end{align*}
and put forth the following proposition.

\begin{proposition}\label{formal refine}
The formal spin 3-manifold invariants $Z_{SO_3}(X,\sigma, 2m-1)$ 
are refinements of the formal 3-manifold invariant $ Z_{SU_2}(X, 4m-2)$
in the sense that 
$$
Z_{SU_2}(X, 4m-2) = \sum_{\sigma}
Z_{SO_3}(X,\sigma, 2m-1).
$$
\end{proposition}

\begin{proof}
We begin with a more down-to-earth computation.  Let $(P,A)$ be an 
$SO_3$-connection such that $w_2(P) = b$.  Then we consider the 
sum
\begin{align*}
\sum_{\sigma} \rttau_{X,\sigma} (D_{\rho A})^{2m-1}  = &
\sum_{\ell} \rttau_{X,\sigma+\ell} (D_{\rho A})^{2m-1} \\
= & \left(  \sum_{\ell} (-1)^{b \smile \ell} \right)
\rttau_{X,\sigma} (D_{\rho A})^{2m-1}.
\end{align*}
The second and third sums are over $\ell \in H^1(X;\zmod2)$.  The
second equality follows from Theorem \ref{properties of q}.  It is
easy to see that the third sum is equal to zero if $b \neq 0$ and is equal to 
$\# H^1(X;\zmod2)$ if $b = 0$.  The upshot of this computation is that,
in summing over the spin structures, the $SO_3$-connections with 
non-trivial $w_2$ contribute nothing and all non-zero contributions
come from the $SO_3$-connections with trivial $w_2$.  In particular,
\begin{align}\label{formal so3}
& \sum_{\sigma} Z_{SO_3}(X,\sigma, 2m-1) = \\
\#( H^1(X;\zmod2) )&  \int_{\overline{\Cat{SO_3,0}{X}}} 
\quad f_{\sigma} \quad \mu_{S0_3,0}(X). \nonumber
\end{align}
Notice that over $\overline{\Cat{SO_3,0}{X}}$ the function $f_{\sigma}$ is 
independent of $\sigma$.  Now if we consider
$$
Z_{SU_2}(X, 4m-2) = 
\int_{\overline{\Cat{SU_2}{X}}} \quad f' \quad \mu_{SU_2}(X)
$$
then \eqref{formal int2}, \eqref{def f'}, and \eqref{formal so3} 
imply the proposition.
\end{proof}

\section{The Quantum Hamiltonian Theory}\label{quantum hamiltonian}

\subsection{Geometric quantization with K\"{a}hler polarization}
\label{2.2.1}

In this section we review the procedure of geometric quantization
when the symplectic manifold has a K\"{a}hler structure.  Let 
$(M, \omega)$ be a smooth complex manifold with a positive K\"{a}hler form 
and let $\Line{}{} \rightarrow M$ be a hermitian line bundle 
with a unitary connection $\nabla$.  Furthermore, we require
that $\nabla \circ \nabla = \Omega = - 2 \pi i \omega$.
To conform with the literature (e.g. \cite{Wo}, \cite{GS}) we
call $\Line{}{}$ the ``pre-quantum line bundle''.  Notice that,
since $\omega$ represents $c_1(\Line{}{})$ in deRham cohomology,
it takes integer values when integrated over smooth 2-cycles.
This is already a ``quantization'' condition of sorts.

The connection $\nabla$ determines a unique holomorphic structure
$\bar{\d}_{\Line{}{}}$ on the pre-quantum line bundle.  This 
follows from the the fact that $\Omega$ is a $(1,1)$ differential 
form in the bi-grading determined by the complex structure on 
$M$.  Thus, one way to obtain a Hilbert space from this system would be
to take the holomorphic sections $\Hilb' = 
H^0(M;\Line{}{})$.  The natural hermitian form for $\Hilb$ is
$$
\pairing{s_1}{s_2}_{\Hilb} = 
\int_M \pairing{s_1}{s_2}_{\Line{}{}} \cdot
\omega^{\dim_{\C}M}  .
$$
However, certain examples (cf. \cite{GS}) tell us that
this is often not the ``right'' Hilbert space.  A slight correction
must be made and we explain this next.

In the most favorable cases, we can include the {\it metaplectic correction} 
in the geometric quantization procedure. 
Let $K = \bigwedge^{\text{top}}(T^{(1,0)}M)^*$ denote the 
canonical line bundle of $M$.  When the real tangent bundle of $M$ has trivial 
2nd Stiefel-Whitney class -- or what is the same, the holomorphic
tangent bundle has 1st Chern class that is trivial under mod-2 reduction --  
there exists a holomorphic line bundle $K^{1/2} \rightarrow M$ such that 
$(K^{1/2})^{\otimes 2} = K$.  In fact, given the existence of a $K^{1/2}$, 
the set of equivalence classes of such line bundles is affine over $H^1(M,\zmod2)$.  
In the literature $K^{1/2}$ is called the ``bundle of half-forms''.  In many cases, 
the ``right'' Hilbert space is $\Hilb = H^0(M;\Line{}{}\otimes K^{1/2})$.

On an affine symplectic vector space, for instance, 
the inclusion of the metaplectic correction is critical to identifying
Hilbert spaces associated to different affine K\"ahler
structures.  More recently, \cite{Hall}, its been discovered that
in certain instances the metaplectic correction is required to identify the 
respective Hilbert spaces associated with a K\"ahler structure and a
{\it real polarization}.   We feel that these cases provide strong 
motivation for including the correction, whenever possible, in our quantization
procedure.

\subsection{Geometric structures on the moduli stack}

Let $Y$ be a closed, 2-manifold with spin structure
$\sigma$ and let $G$ be a compact Lie group.  For simplicity
we further require that $G$ be connected.  The moduli stack $\Mod{G}{Y}$
is, of course, independent any Riemannian or spin structure on $Y$
but given such structures we obtain certain geometric structures on the moduli 
space.  We discuss this next.     
  
Fix a Riemannian structure on $Y$ and let $Q \rightarrow Y$ be a principal 
$G$ bundle.  Any Riemannian 2-manifold has a unique complex structure
determined by it Hodge star operator.  It is trivially integrable.  According to a 
classical theorem of Narasimhan and Seshadri \cite{NS} (see also \cite{AB}), 
a complex structure on $Y$ induces a complex structure on the moduli space 
$\Mod{}{Q}$ of flat connections in the following way.  We let $G_{\C}$ denote the 
complexification of $G$ and we let $Q_{\C}$ denote the natural extension of 
$Q$ to a $G_{\C}$ bundle.  There is an identification between 
$\Mod{}{Q}$ and the moduli space $\Mod{}{Q_{\C}} = \Cat{}{Q_{\C}}^{ss} 
/ \Mor{}{Q_{\C}}$ of semi-stable holomorphic structures on $Q_{\C}$.  The 
space $\Mod{}{Q_{\C}}$ has a natural complex structure.  Indeed, over a 
smooth point $B$,  the holomorphic tangent space is modeled on 
$H^1(Y; \bar{\d}_B)$, the first cohomology group of the complex
\begin{equation}\label{cmplx}
0 \rightarrow \Omega^0 (Y; \text{ad}Q_{\C})
\xrightarrow{\bar{\d}_B}
\Omega^{0,1}(Y; \text{ad}Q_{\C}) \rightarrow 0 .
\end{equation}  
Here the superscript $(0,1)$ implies that these are the anti-holomorphic
1-forms with respect to the complex structure on $Y$ and $\bar{\d}_B$
is the unique holomorphic structure on $Q_{\C}$ determined by the
flat connection $B$.     

Recall from \ref{Classical hamiltonian} that, given a real representation 
$\rho$ of rank zero, there is a hermitian line bundle $\Line{\rho}{(Q)} \rightarrow
\Mod{}{Q}$ with a unitary connection $\nabla$.  If $\Omega$ is the 
curvature of $\nabla$ then over smooth points of the moduli space,
$\omega = i \Omega / 2 \pi$ defines a symplectic structure.  With respect 
to the complex structure on $\Mod{}{Q}$ the symplectic
form $\omega$ has bigrading $(1,1)$.  Thus $\Line{\rho}{(Q)}$ 
is a holomorphic line bundle and we are in the scenario described in the 
previous section.  

As mentioned above, to obtain the ``right'' Hilbert space we require the 
canonical line bundle $K \rightarrow \Mod{}{Q}$ and then a square root
of that bundle.  We discuss this next.  For the rest of this section we assume 
the genus of $Y$ is greater than one.

The zeroth cohomology $H^0(Y;\bar{\d}_B)$ of the complex 
\ref{cmplx} is the Lie algebra of the stabilizing subgroup of 
$\bar{\d}_B$.  If $B$ represents a smooth point then the stabalizing 
subgroup consists of the center of $G$ so that $H^0(Y;\bar{\d}_B)
= 0 $.  Thus, at a smooth point $B$ we have the equivalence of lines
\begin{equation}\label{eqv lines}
\Det^{-1}_{\bar{\d}_B} = \Det H^0(Y;\bar{\d}_B) \otimes
\Det^{-1} H^1(Y;\bar{\d}_B) =
\Det^{-1} H^1(Y;\bar{\d}_B) .
\end{equation}
Given that $H^1(Y;\bar{\d}_B) = T^{1,0} \Mod{}{Q}$ we see that
the far right hand side of \eqref{eqv lines}is the fiber of the canonical 
bundle $K$ at $B$.  If we let $\Det^{-1}_{Q_{\C}} \rightarrow 
\Cat{}{Q_{\C}}$ denote the inverse determinant line bundle whose fiber 
at $\bar{\d}_B$ is the left hand side of \eqref{eqv lines}, then 
$\Det^{-1}_{Q_{\C}}$ descends to $\Mod{}{Q}$.  Over smooth points 
the descendant line bundle $\Det^{-1}_{Q_{\C}} \rightarrow \Mod{}{Q}$ 
is equivalent to the canonical bundle $K \rightarrow \Mod{}{Q}$.  The 
inverse determinant line bundle definition of $\Det^{-1}_{Q_{\C}}$ is 
rather convenient in the setting of the Hamiltonian spin Chern-Simons theory 
with its Pfaffian line bundles.  For this reason with work with $\Det^{-1}_{Q_{\C}}$ 
instead of $K$. 

As a guide to finding a square root of $\Det^{-1}_{Q_{\C}}$ we point out that we 
may write
\begin{equation}\label{dbar wrt dirac}
\bar{\d} = D \otimes K_Y^{-1/2}
\end{equation}
where the left hand side is the usual Dolbeault operator on $Y$ and the
right hand side is the Dirac operator twisted by the square root
of the canonical line bundle.  Actually, $K_Y^{-1/2}$ is the spinor bundle 
$S_Y^-$ transposed to the holomorphic setting.
The same follows when the two operators are twisted by any vector bundle
with connection; in particular it  follows when they are twisted by 
$\text{ad}Q_{\C}$.  We offer the following proposition.  For simplicity
we now assume $G$ is connected.

\begin{proposition}\label{rtK}
Let $\rho_{\ad} = \ad - \dim G$ be
the adjoint representation of $G$ minus the trivial representation
of rank $\dim G$.  Then  there exists a connection preserving isometry 
between the two inverse determinant lines $\Det^{-1}_{Q_{\C}} 
\rightarrow \Mod{}{Q}$ and $(\Line{\rho_{\ad}}{(Q)})^{\otimes 2} 
\rightarrow \Mod{}{Q}$.  
\end{proposition}

\begin{remark}
As will be clear from the proof, the above proposition applies more
generally.  Indeed, $\rho_{\ad} Q$ may be replaced with any real oriented
vector bundle and $K_Y^{-1/2}$ may be replaced with any fixed complex
line bundle on $Y$.  The upshot of this theorem is that we may model
the canonical line bundle of $\Mod{}{Q}$ on the line bundle
$(\Line{\rho_{\ad}}{(Q)})^{\otimes 2}$ and so we may model the 
bundle of half-forms on the line bundle $(\Line{\rho_{\ad}}{(Q)})$.
Notice that a choice of square root for the canonical bundle $K_Y$
determines a choice of square root for the canonical bundle $K$ of
the moduli space.
\end{remark}

\begin{proof}
As the statement of the proposition suggests, we only prove existence
and do not construct an isomorphism. To prove existence we show that 
the line bundle $\Line{}{} = \Det^{-1}_{Q_{\C}} \otimes 
(\Line{\rho_{\ad}}{(Q)})^{-\otimes 2}$ has
trivial holonomy.  Then there exists a covariantly constant
unitary section which is unique up to a factor in $\T$.

Now we show that the holonomy around any closed path 
in $\Mod{}{Q}$ is the identity.  From \label{Holonomy} we know that the
holonomy around any loop (in the case of inverse determinant
lines) is given by the $\tau$-invariant of some twisted Dirac operator
over $S^1 \times Y$ where the bounding spin structure is placed on $S^1$.  
What we show now is that 
\begin{equation}\label{triv tau}
\tau_{S^1 \times Y}
(D_E \otimes (K_Y^{-1/2} - 1)) = 1
\end{equation}
for any real oriented virtual vector bundle $E\rightarrow S^1 \times Y$ 
with orthogonal connection $\nabla^E$. 
We assume that the metric $S^1 \times Y$ 
is product and furthermore that the metric on $S^1$ is flat.  This caveat is 
part of the hypothesis of the theorem which expresses the holonomy in terms 
of (adiabatic limits of) $\tau$-invariants.

We first claim that the right hand side of \eqref{triv tau} 
in independent of  $\nabla^E$.  Indeed, the formula for the
differential (of the log) of $\tau$ is given by the Atiyah-Patodi-Singer
index theorem \cite{APS1}, and a 
straightforward computation of the local index 
shows that the differential is zero.  From
this we see that the right hand side of \eqref{triv tau} only
depends on the topological type of $E$.  Because we are working over a 
3-manifold, the topological type of $E$ is completely determined by $w_2(E)$. 
(See the appendix of \cite{Je1} for a proof).  The Kunneth formula implies 
$$
H^2(S^1 \times Y;\zmod2) \cong H^2(Y;\zmod2) 
\oplus H^1(S^1;\zmod2) \otimes_{\zmod2} H^1(Y;\zmod2). 
$$
Since both $w_2$ and $\tau$ depend linearly on $E$ we can assume,
without loss of generality, that $w_2(E)$ lies in one of the above summands.
Then we can show that, in each case, the right hand side of \eqref{triv tau} is 1.
First assume that $w_2(E) \in H^2(Y;\zmod2)$.  In this case $E$ extends to 
a bundle $E' \rightarrow D^2 \times Y$ so that the APS index theorem
implies
\begin{align*}
& \tau_{S^1 \times Y}
(D_E \otimes (K_Y^{-1/2} - 1)) = \\
& \exp 2 \pi i \left[ \int_{D^2 \times Y} 
\widehat{A}(\Omega^{D^2 \times Y}) ch(\Omega^{E'})
ch(K^{-1/2} - 1)  \right]_{(0)}
\end{align*}
and another straightforward computation shows that the integral 
is zero.   Now let $\ell$ denote the non-trival element of $H^1(S^1;\zmod2)$
and assume that $w_2(E) = \ell \smile a$ for some 
$a \in H^1(Y;\zmod2)$.  As we did in Section 1.1, identify $a$ and 
$\ell$ with the flat real line bundles that represent them.  Then we
can assume that $E = a \oplus \ell \oplus \ell \otimes a$ so that 
\begin{align*}
& \tau_{S^1 \times Y}(D_E \otimes (K_Y^{-1/2} - 1)) = \\
& \tau_{S^1 \times Y}(D_a \otimes (K_Y^{-1/2} - 1)) \cdot
\tau_{S^1 \times Y}(D_\ell \otimes (K_Y^{-1/2} - 1)) \cdot
\tau_{S^1 \times Y}
(D_{\ell \otimes a } \otimes (K_Y^{-1/2} - 1)).
\end{align*}
The first factor involves only vector bundles that can be extended over 
$D^2 \times Y$ and can be shown to be equal to 1 using the same 
argument as in the case $w_2(E) \in H^2(Y;\zmod2)$.  To deal with the
last two factors we note that twisting by $\ell$ is equivalent to placing
the non-bounding spin structure on $S^1$, so that
$$
\tau_{S^1 \times Y}(D_E \otimes (K_Y^{-1/2} - 1)) =
(-1)^{\text{index}D \otimes (K_Y^{-1/2} - 1)} 
(-1)^{\text{index}D_a \otimes (K_Y^{-1/2} - 1)}. 
$$ 
Of course the index cannot detect the factor $a$ so that the two
indices above are equal and so cancel each other out (mod 2).  Finally
we see that the right hand side of \eqref{triv tau} is 1 in this
case, as well.  
\end{proof}

\subsection{The moduli spaces for $G = SU_2$, $G = SO_3$ }

We review some salient features of the moduli spaces 
$\Mod{SU_2}{Y}$ and $\Mod{SO_3}{Y}$.  In particular, we
restrict ourselves to the case of connected oriented 2-manifolds of genus 2 
or higher.

For any connected group $G$, the topological equivalence classes of principal
$G$-bundles over $Y$ are parametrized by the abelian group $H^2(Y; \pi_1 G) 
\cong \pi_1 G$.  Indeed, if we remove a disc $D$ from $Y$, then $Y \setminus
D$ is homotopic to its 1-skeleton over which any $G$-bundle is trivial.  Thus,
up to isomorphism, any $G$-bundle on $Y$ is determined by the homotopy
class of the clutching map $g : \d D \rightarrow G$ that glues the
trivial bundle over $Y \setminus D$ to the trivial bundle over $D$.  The 
elements of $\pi_1 G$ are in one-one correspondence with the connected
components of $\Mod{G}{Y}$ since each topological type determines
a component of the moduli stack.  In particular, $\Mod{SU_2}{Y}$ has one 
component while $\Mod{SO_3}{Y}$ has two components: one for the each
of the two possible 2nd Stiefel-Whitney classes.  We denote these components
by $\Mod{SO_3, w_2}{Y}$ for $w_2 = 0,1$.

For any compact Lie group $G$ each connected compontent of the moduli space 
decomposes into a disjoint union of strata.  The stratum of highest dimension is 
a smooth manifold and is a dense open subset of the component with respect to 
to quotient topology.  For $SU_2$ the top stratum consists of the 
{\it irreducible} flat connections.  To explain what we mean by ``irreducible''
we remind the reader that a flat $G$-connection on $Y$ determines a 
representation $\pi_1 Y \rightarrow G$.  It is this representation that
is irreducible.  The other strata of the $SU_2$ moduli space consist of 
{\it reducible} flat connections, in particular those that reduce to flat $\T$
and $\zmod2$ connections.  Recall that any maximal torus of $SU_2$ is
isomorphic to $\T$ and the center of $SU_2$ is isomorphic to $\zmod2$;
this explains the appearance of $SU_2$-connections that reduce to these 
Lie groups.   

For $\Mod{SO_3,0}{Y}$ and $\Mod{SO_3,1}{Y}$ moduli spaces the top 
strata also consist of irreducible flat connections.  The other strata consist 
of flat connections that reduce to flat
$O_2$, $SO_2$, $\zmod2$ and $\zmod2 \times \zmod2$ connections.
Recall that $SO_2$ embeds into $SO_3$ as rotations about a fixed 
axis and $O_2$ embeds as rotations about a fixed axis as well as 
such rotations composed with a $180^{\circ}$ rotation about some
perpendicular axis.  The group $\zmod2$ embeds as a subgroup of 
$SO_2$ and $\zmod2 \times \zmod2$ embeds as compositions
of $180^{\circ}$ rotations about three mutually orthogonal axes.
This explains the appearance of $SO_3$-connections that reduce to
these Lie groups.  A trivial $SO_3$ bundle supports connections that
reduce to each of these subgroups, however a non-trivial $SO_3$ bundle
does not support connections that reduce to flat $SO_2$ or $\zmod 2$
connections.

The complex structure of the $SU_2$ moduli space comes by way of
extending the $SU_2$ bundles to $(SU_2)_{\C} = SL_2(\C)$ bundles 
and then identifying the moduli space of flat connections with the 
moduli space of semi-stable holomorphic structures.  Similarly, we extend
the $SO_3$ bundles to $(SO_3)_{\C} = \P L_2(\C)$ bundles and 
use the identification of flat connections with holomorphic structures.

\subsection{Flat $SO_3$ and Yang-Mills $U_2$ connections}

We now address the apparent equality of dimensions for the top strata of the
$SU_2$ and $SO_3$ moduli spaces.  To do so we have to extend our
consideration to {\it Yang-Mills} connections on principal $U_2$
bundles over $Y$ of which flat $SU_2$ connections are a subset.
These connections were the main object of study in \cite{AB} and much
of what we have to say about them comes from that reference.

Consider the short exact sequence of compact Lie groups
\begin{equation}
1 \rightarrow \T \longrightarrow U_2 
\longrightarrow \P U_2 \rightarrow 1.
\end{equation}
The first map is multiplication by the $2 \times 2$ identity matrix and maps
$\T$ into the center $C \subset U_2$.  The second map is projectivization.  To 
connect to $SO_3$ we point out that $\P U_2 \cong SO_3$.  Indeed, the adjoint 
action of $U_2$ on the Lie algebra preserves the decomposition $\Lie{u}_2 = 
\Lie{c} \oplus \Lie{su}_2$; it acts trivially on $\Lie{c}$ and the non-trivial 
adjoint action on $\Lie{su}_2$ factors through $\P U_2$.  A choice of 
identification $\Lie{su}_2 \rightarrow R^3$ induces a bijection $\P U_2 
\rightarrow SO_3$.     

If $A$ is a $G$-connection over the oriented Riemannian 2-manifold $Y$,
then it is {\it Yang-Mills} if $d_A * \Omega^A = 0$, where $*$ denotes
the Hodge star operator of the metric on $Y$.  This implies
that $\Tr (i \Omega^A / 2\pi )$ is the unique harmonic form that represents 
the first Chern class in deRham cohomology.  The gauge group $\Mor{G}{Y}$
acts on the space of Yang-Mills connections and we can therefore talk about the
corresponding moduli space.
The moduli space of Yang-Mills $U_2$ connections is a well-studied space,
and we use the identification $\P U_2 = SO_3$ to our advantage in 
quantizing our spin-Chern-Simons theory.  We explain how next.  

Consider the short exact sequence of compact Lie groups
\begin{equation}\label{ses}
1 \rightarrow \zmod2 \longrightarrow U_2 
\longrightarrow \T \times \P U_2 \rightarrow 1.
\end{equation}
The first map is multiplication of $\pm 1$ by the $2 \times 2$ identity matrix 
and the second map is the cartesian product of the maps $\det : U_2 
\rightarrow \T$ and $\P : U_2 \rightarrow \P U_2$.  From this short exact 
sequence we see that a Yang-Mills $U_2$ connection $(P,A)$ is determined by 
the Yang-Mills connection $(\Det P , \Det A)$ and the flat connection 
$(\P P, \P A)$ up to a $\zmod2$-connection.  To be precise, the category 
$\Cat{\zmod2}{Y}$ of $\zmod2$-connections acts on $\Cat{U_2}{Y}$
by ``tensor product''.  Indeed, for $Z \in \Cat{\zmod2}{Y}$ and 
$(P,A) \in \Cat{U_2}{Y}$ we obtain a new $U_2$-connection:
$$
(Z \otimes P, Z \otimes A) = (Z \times_Y (P,A))/ \zmod2
$$     
where the quotient is taken with respect to the diagonal $\zmod2$ action
on the fiber product.  This is the principal bundle analogue of tensoring
a flat real orthogonal line bundle with a hermitian rank 2 vector bundle,
which justifies our ``tensor product'' notation.  We will use this notation
throughout the rest of this section.  The important fact to observe is that
tensoring $A$ with a $\zmod2$-connection preserves $\Det A$ and 
$\P A$.  This is what we mean when we say that $\Det A$ and $\P A$
determine $A$ up to a $\zmod2$-connection. 

If we fix a Yang-Mills $\T$-connection $(T, a)$ we can consider the category 
$\Cat{U_2}{T,a}$ of $U_2$ connections $(P,A)$ such that $(\Det P , 
\Det A) = (T,a)$.  The morphisms $\Mor{U_2}{T,a}$ of this category 
are those elements $\phi$ of $\Mor{U_2}{Y}$ such that $\Det \phi =
\text{id}_T$; and so we may consider the moduli space $\Mod{U_2}{T,a}
= \Cat{U_2}{T,a} / \Mor{U_2}{T,a}$.  Now an element $A \in 
\Cat{U_2}{T,a}$ is determined by $\P A$ up to a $\zmod2$-connection.
We will use this fact to identify the moduli space $\Mod{SO_3,w_2}{Y}$
as the quotient of some $\Mod{U_2}{T,a}$ with respect to $\Mod{\zmod2}{Y}$.  
To discern the correct choice of $(T,a)$ we point out that for a $U_2$ principal
bundle $P \rightarrow Y$, $w_2(\P P) = c_1(P)$ (mod 2) and that 
$c_1(P) = c_1(\Det P)$.  Thus we see that if $(T,a)$ has even (resp. odd)
degree then $\Mod{SO_3,0}{Y}$ (resp. $\Mod{SO_3,1}{Y}$) is identified with 
$\Mod{U_2}{T,a}/ \Mod{\zmod2}{Y}$.  In particular, if we fix $(T,a)$ to be the 
trivial line bundle with trivial connection, then it is clear that 
$\Mod{U_2}{T,a} = \Mod{SU_2}{Y}$ and we can identify the quotient of the
latter with $\Mod{SO_3,0}{Y}$ as well.
      
\subsection{The lift of the $\Mod{\zmod2}Y$ action}

Let $Y$ be a closed, oriented, genus-$g$ 2-manifold with fixed metric and 
spin structure $\sigma$, and let $\beta : U_2 \rightarrow SO_3$ be the 
projective homomorphism composed with some isomorphism $\P U_2 
\rightarrow SO_3$.  We also fix a real, rank zero virtual representation 
$\rho $ of $SO_3$ and a degree-1 Yang-Mills $\T$-connection $(T,a)$. 
Throughout this section we use the following shorthand.  We let $\Cat{d}{Y}$  
stand in for either $\Cat{SU_2}{Y}$ or $\Cat{U_2}{T,a}$;  in particular 
$d=0$ for the category of $SU_2$ connections and $d=1$ 
for the category of $U_2$ connections whose determinants are equal to $(T,a)$.  
We let $\Cat{w_2}{Y}$ stand in for $\Cat{SO_3, w_2}{Y}$.  We use a corresponding 
shorthand for the $SO_3$ and $U_2$ moduli stacks.  To avoid typographical 
redundancy and emphasize its dependence of the spin structure, we drop the
 ``$Y$'' from the notation for the Pfaffian lines, writing $\Line{(\rho)}
{\sigma}$ for $\Line{\rho}{(Y)}$ and $\Line{\rho \circ \beta}{}$ for 
$\Line{\rho \circ \beta}{(Y)}$.   

In this section we take advantage of known results for the vector spaces 
$H^0(\Mod{d}{Y}; \Line{\rho \circ \beta})$, $d=0,1$.  In particular, 
the dimensions of these spaces are the well-known Verlinde formulae; the
vector space for $d = 1$ (and its corresponding Verlinde formula) is often 
referred to as the {\it twisted} case.  We will use these previously
obtained Verlinde formulae to obtain formulae for the dimensions
of $H^0(\Mod{w_2}{Y}; \Line{\rho}{(\sigma)})$ for $w_2=d$.  
The formulas we obtain are not new and correspond to the formulas in 
\cite{BM} and \cite{AM}.  The novelty here is our approach to the 
computation.  Our vector spaces are derived from line bundles over the 
$SO_3$ moduli stacks while the \cite{AM} computations are done in the 
context of $U_2$ moduli stacks.  In \cite{AM} the authors consider
a central extension of $\Mod{\zmod2}{Y}$ and its action on  
$H^0(\Mod{d}{Y}; \Line{\rho \circ \beta})$.  They derive the 
vector spaces as the irreducible components of this action.  See Section 
3.1 or the cited reference for more details.
 
Our main concern here is the action of $\Mod{\zmod2}Y = H^1(Y;\zmod2)$  
on $\Mod{d}Y$.  This action really follows from an action of 
$\Cat{\zmod2}Y$ on $\Cat{d}Y$ given by the tensor product action 
described above:
\begin{align*}
\Cat{\zmod2}Y \times \Cat{d}Y & 
\longrightarrow \Cat{d}Y \\
( Z , (P,A) ) &
\longmapsto (Z \otimes P, Z \otimes A)
\end{align*}  
This categorical action descends to an action on the moduli stacks
\begin{align*}
\Mod{\zmod2}Y \times \Mod{d}Y & 
\longrightarrow \Mod{d}Y \\
( [Z ] , [P,A] ) &
\longmapsto [(Z \otimes P, Z \otimes A)]
\end{align*}
In what follows we will make use of the identification
$$
\Mod{d}Y / \Mod{\zmod2}Y \cong \Mod{w_2}Y .
$$
Before we do anything more we characterize the fixed points of $[Z]$ 
in $\Mod{d}Y$ by their images in $\Mod{w_2}Y$.  
In particular we have the following proposition.

\begin{proposition}[Characterization of the Fixed Points]\label{fxd pts}
The equivalence class $[P,A]$ is a fixed point of $[Z]$ if and only 
if $(\beta P, \beta A)$ can be reduced to an $O_2$-connection
$(Q,B)$ such that $\Det Q \cong Z$.
\end{proposition}

Here $\Det Q$ denotes the orientation bundle of $Q$ with its natural flat 
connection.  This characterization is extremely useful in the forthcoming
compuation of $\dim H^0(\Mod{w_2}{Y}; \Line{\rho}{(\sigma)})$. 

\begin{proof}
We first prove sufficiency.  Let $Q \subset \beta P$ be an $O_2$ subbundle
such that $\beta A$ preserves $Q$ under parallel transport.  We need to 
construct a morphism of pairs 
$$
\Phi: (\Det Q \otimes P, \Det Q \otimes A) \longrightarrow
(P,A) .
$$
Let $R$ denote the rotation
$$
R =
\begin{pmatrix}
-1 & 0 & 0 \\
0 & -1 & 0 \\
0 & 0 & 1
\end{pmatrix}
\in O_2 \subset SO_3
$$
and let $\tilde{R} \in SU_2$ be chosen so that 
$\beta (\tilde{R}) = R$.  If $q, \beta p$ lie in the same fiber
$Q_x$ over $x \in Y$ then we set
$$
\Phi( \Det q \otimes p) = p \cdot \tilde{R} \cdot \frac{\Det q}{\Det \beta p}
$$   
where 
$$
\frac{\Det q}{\Det \beta p} \in \{ \pm 1 \}
$$ 
is the translation
that carries $\Det \beta p$ to $\Det q$.  We let the reader check that
this is a well-defined morphism over for all such $q, p$.  Equivariance now 
defines $\Phi$ at all other elements of $\Det Q \otimes P$.  Clearly
$\Phi$ descends to an automorphism of $(\beta P, \beta A)$ so that
it must send $\Det Q \otimes A \mapsto A$.  This proves sufficiency.

To prove necessity we assume that there exists a morphism
$$
\Phi : (Z \otimes P, Z \otimes A) \longrightarrow (P,A)
$$
If the stabilizing subgroup of $(P,A)$ is isomorphic to 
$U_1$ or $\zmod2$ we can always chose $\Phi$ so that  
$$
(id_Z \otimes \Phi) \circ \Phi = - id_P.
$$
If the stabilizing subgroup is isomorphic to $U_2$ then 
their is no non-trivial $Z$ that stabilizes $(P,A)$.
Thus, in the non-trivial case, we can always chose $\Phi$ 
so that it descends, via $\beta$, to an automorphism
$$
\phi : (\beta P, \beta A) \longrightarrow (\beta P, \beta A)
$$
such that $\phi \circ \phi = id_{\beta P}$.  If $\phi = 
id_{\beta P}$ then its clear that $Z$ is trivializable.  Assuming $\phi 
\neq id_{\beta P}$ it determines an $O_2$ sub-bundle 
$$
Q = \{ q \in \beta P | \phi (q) = q \cdot R \}
$$
which inherits a connection $B$ from $\beta A$ (since $\phi$ preserves
$\beta A$).  Our final task is to construct an isomorphism between $Z$ and 
$\Det Q$.  From the sufficiency argument we know that $\phi$ can be 
lifted to a morphism
$$
\Phi' : (\Det Q \otimes P, \Det Q \otimes A) \longrightarrow
(P,A)
$$ 
For $\Det q \in \Det Q_x$ and $z \in Z_x$we define the map
$$
\Det q \mapsto z \cdot \frac{\Phi'(\Det q \otimes p)}{\Phi(z \otimes p)}
$$
where 
$$
\frac{\Phi'(\Det q \otimes p)}{\Phi(z \otimes p)} \in \{ \pm 1 \} \subset U_2
$$
is the translation that takes $\Phi'(\Det q \otimes p) \in P_x$ to 
$\Phi(z \otimes p) \in P_x$.  We let the reader check that the map
is well-defined despite the choice of $z$.  This completes the proof.
\end{proof}

We would like to see if and then how this lifts to the line bundles 
$\Line{\rho \circ \beta}$.  As before, to see the action we
start on the categorical level.  Indeed the natural isometry
\begin{equation}
\Line{\rho \circ \beta}{(P,A)}
\longrightarrow 
\Line{\rho \circ \beta}{(Z \otimes P, Z \otimes A)} 
\end{equation}
is almost trivial thanks to the natural identification
\begin{align*}
\beta (Z \otimes P) & \longrightarrow \beta (P) \\
\beta(z \otimes p) & \longmapsto \beta (p) 
\end{align*}
This naturally gives us the commutative diagram
\begin{equation}\label{comm diag1}
\begin{CD}
\Line{\rho \circ \beta}{}  @>>>  \Line{\rho}{} \\
@VVV                               @VVV \\
\Cat{d}{Y}  @> \P >>   \Cat{w_2}{Y}.
\end{CD}
\end{equation}

\begin{remark}
This lift has a couple of nice properties.  First of all, the lift induces a ``group'' 
action of $\Cat{\zmod2}Y$ on $\Line{\rho \circ \beta }$.
Thus $\Line{\rho \circ \beta }$ is a $\Cat{\zmod2}Y$
equivariant line-bundle.  Indeed, this is what is implied
by \eqref{comm diag1}.
Second of all, the lift is covariant with respect to the
natural connection on the bundles and commutes with 
the action of the bundle morphisms (after natural 
identifications).  Thus the group action of $\Cat{\zmod2}Y$
not only lifts to the Pfaffian line bundle, it does so
covariantly with respect to the connection over the moduli stack.
Thus have the following commutative diagram:
\begin{equation}\label{comm diag2}
\begin{CD}
\Line{\rho \circ \beta}{}  @>>>  \Line{\rho}{} \\
@VVV                               @VVV \\
\Mod{d}{Y}  @> \P >>   \Mod{w_2}{Y}.
\end{CD}
\end{equation}
\end{remark}
Now we need only compute the actions over fixed points
to know exactly how the actions lift over the whole bundle.
With that in mind, we wish to see how this action lifts 
when we have a fixed point $[Z \otimes P ,Z \otimes A ] 
= [P,A]$ on  $\Mod{d}Y$.  In this case there must be a 
morphism
$$
\Phi : (Z \otimes P , Z \otimes A) 
\longrightarrow (P,A)
$$   
which, under the natural identification 
$\beta (Z \otimes P) = \beta P$ descends to an automorphism
$$
\phi : (\beta (P), \beta (A) ) \longmapsto
(\beta (P), \beta (A) ).
$$
Notice that, under the natural identification
$Z \otimes Z \otimes P = P$
we can compose $\Phi$ with the morphism
$$
id_Z \otimes \Phi : 
(Z \otimes Z \otimes P, Z \otimes Z \otimes A)
\longmapsto (Z \otimes P , Z \otimes A)
$$
so that we get an automorphism
$$
( id_Z \otimes \Phi ) \circ \Phi: (P,A) 
\longmapsto (P,A)
$$
which, via $\beta$, descends to the automorphism 
$\phi \circ \phi$.  

We are finally in a good position to compute the action 
of $[Z] \in \Mod{\zmod2}Y$ on a fixed point.  
The result is the following

\begin{proposition}[The Lift to the Pfaffian Line]
The lift of $[Z]$ to the line $\Line{\rho \circ \beta}{[A]}$ 
over a fixed point $[A]$ is multiplication by 
$$
q (\sigma, \rho (\beta P), Z) = 
(-1)^{w_2 (\rho \circ \beta P) + w_2(\rho)(\text{ind}_2(D_{\sigma + Z}) 
- \text{ind}_2(D_{\sigma }) ) }
$$
where $w_2 (\rho \circ \beta P) \in \zmod2$ denotes the
invariant $w_2( \rho \circ \beta P) \frown [Y]$ and $w_2(\rho) \in \zmod2$
is obtained from the map $w_2 : RO(SO_3) \rightarrow H^2(BSO_3;\zmod2) 
\cong \zmod2$.
\end{proposition}

\begin{proof}
By the proposition \ref{fxd pts} we have that
there is a sub-bundle $(Q,B) \subset (\beta P, \beta A)$ 
such that $\Det Q \cong Z$.  We also have that $\phi$ -- 
the descendent (via $\beta$) of the morphism 
$$
\Phi: (Z \otimes P, Z \otimes A) \longrightarrow (P,A)
$$
-- restricts to $(-id_Q)$ on $Q$.  
The action on the line is multiplication by the $\rttau$-invariant of
$$
(Q, B) \times_{(-id_Q) } S^1_b ,
$$   
i.e. the bundle with connection gotten gluing the ends of $(Q,B) \times
[0,1]$ together with $\phi|Q = -id_Q$.  Notice that in doing so
we give $S^1$ the {\it bounding} spin structure.  Recall 
that this is what is required to compute the {\it trace} (as opposed to the 
{\it super-trace}).  Then, letting 
$$
i: O_2 \hookrightarrow SO_3
$$
denote the standard inclusion homomorphism, we want to compute
$$
\rttau ( \rho \circ i (Q,B) \times_{(-id_Q)} S^1_b)   = 
\rttau ( \rho \circ i (Q,B) \times_{id_Q} S^1_{nb})   = 
(-1)^{\text{\rm ind}_2( D_{\sigma} \otimes \rho \circ i Q )} 
$$
Note that the mod-2 index is independent of $B$ (as it is a topological
invariant).  All that remains is to compute 
\begin{align}
\text{\rm ind}_2(D_{\sigma} \otimes  \rho \circ i (Q) )& = 
w_2 \rho \circ i (Q)  + w_2(\rho)(\text{\rm ind}_2(D_{\sigma + \Det Q}) 
- \text{\rm ind}_2(D_{\sigma }) ) \\
 & = w_2 \rho (\beta P ) + w_2(\rho)(\text{\rm ind}_2(D_{\sigma + Z}) 
- \text{\rm ind}_2(D_{\sigma }) ) 
\end{align}
where the first equality follows from the $KO$-theoretic decomposition
of $\rho \circ i (Q)$ and the second equality follows from $\Det Q \cong Z$.
\end{proof}

Of course, this is not the only way in which one can lift the action of 
$[Z]$ to $\Line{\rho \circ \beta}$.  Indeed,
we can also take the (rather natural) lift described above and multiply
it by the scalar factor $q (\sigma ,\rho \circ \beta (P), Z)$.  
This particular action, according to the previous theorem, lifts to a trivial 
action over the fixed points.  The first action, which is more natural in the
context of our spin Chern-Simons field theory, we denote by 
$[Z]_{CS, \sigma}$ while the second we denote by $[Z]_B$ in honor of 
A. Beauville who computed the trace of this particular action on the vector
space $H^0(\Mod{d}{Y}; \Line{\rho \circ \beta} )$ \cite{Be}.  
This trace was also independently computed by J. Andersen and G. 
Masbaum \cite{AM} and T. Pantev \cite{Pa} in the non-trivial
$SO_3$-bundle case.

\begin{theorem}[A. Beauville]
$\Tr [Z]_B =  ( \lambda (\rho ) + 1 )^{g - 1}$
where $\lambda (\rho) \in \Z$ is obtained from the map 
$\lambda : RO(SO_3) \rightarrow E^4(BSO_3) \cong \Z$.
\end{theorem}
This computation makes use of the Lefschetz-Riemann-Roch fixed point
formula \cite{AS} and algebro-geometric results characterizing the fixed points
as certain abelian varieties \cite{NR}.  In the case $w_2 \neq 0$ these are the 
{\it Prym varieties}, and in the case $w_2 = 0$ these are the {\it Kummer
varieties}.  The singularity of the fixed points in the case $w_2 = 0 $
is dealt with by ``transferring'' the computation to the $w_2 \neq 0$
moduli space via the {\it Hecke correspondence}.  For details we refer the
reader to Beauville's paper \cite{Be}.

An easy corollary is then

\begin{corollary}
$\Tr [Z]_{CS, \sigma} = 
q (\sigma , \rho (\beta P), Z ) ( \lambda (\rho ) + 1 )^{g -1}$  
\end{corollary}

We are finally in a good position to compute the dimension of 
$H^0(\Mod{w_2}{Y}; \Line{\rho}{(\sigma)})$.  Indeed, the commutative 
diagram \eqref{comm diag2} tells us that there is a one-one correspondence 
between $H^0(\Mod{w_2}{Y};  \Line{\rho}{(\sigma)} )$ and the 
$\Mod{\zmod2}Y$ invariant subspace of $H^0(\Mod{d}{Y};
\Line{\rho \circ \beta} )$.  The dimension of the latter is just the trace of 
the projection
$$
P_{\sigma} = \frac1{2^{2g}} \sum_{ [Z] } [Z]_{CS, \sigma}
$$
where, unless there is notation to indicate otherwise, the sum is over
all $[Z] \in \Mod{\zmod2}Y$.  Based on this we have

\begin{proposition}\label{dim form}
If $Y$ is a genus-$g$ 2-manifold then 
\begin{align*}
& \dim H^0(\Mod{w_2}{Y}; \Line{\rho}{(\sigma)} ) = \\
\frac1{2^{2g}} 
( \dim H^0(\Mod{d}{Y}; \Line{\rho \circ \beta} ) + 
& (-1)^{w_2 (\rho \circ \beta P)  } 
( (-1)^{\epsilon (\sigma )} 2^g -1) (\lambda (\rho) + 1)^{g-1})
\end{align*}
where $\epsilon (\sigma )$ is the Arf-invariant of $\sigma$.
\end{proposition} 
 
 \begin{remark}
 The Arf-invariant has an index-theoretic formulation that is perhaps
 more appropriate considering the bent of our approach.  Indeed, 
 $\epsilon (\sigma ) = \text{\rm ind}_2 (D_{\sigma } )$; that is the Arf-invariant
 is just the mod-2 index of the (uncoupled) Dirac operator associated
 to $\sigma$.  For future reference we recall the well known fact that
 on a genus-$g$ surface there are $(2^{2g-1} + 2^{g-1})$ spin structures
 $\sigma $ for which $\epsilon (\sigma ) = 0$ and $(2^{2g-1} - 2^{g-1})$
 spin structures $\sigma'$ for which $\epsilon (\sigma' ) = 1$. 
 \end{remark}
 
\begin{proof}
From the discussion above we see that the dimension of 
$H^0(\Mod{w_2}{Y}; \Line{\rho} )$ is given by the formula
$$
\Tr (\frac1{2^{2g}} \sum_{ [Z]} [Z]_{CS, \sigma} )
 = 
\frac1{2^{2g}} \sum_{ [Z] } \Tr( [Z]_{CS, \sigma} )
$$
For $[Z] \neq 0$ we have that 
$$
\Tr [Z]_{CS, \sigma} = 
q (\sigma, \rho (\beta P), Z ) ( \lambda (\rho ) + 1)^{g-1}
$$
and, of course, if $[Z] = 0$ we have 
$$
\Tr [Z]_{CS, \sigma} = \text{dim}
H^0( \Mod{d}{Y}; \Line{\rho \circ \beta }{} ) .
$$
Combining this with the formula above we have that the dimension
of $H^0(\Mod{w_2}{Y}; \Line{\rho} )$ is given by
\begin{equation}\label{dim1}
\frac1{2^{2g}}
\left( \text{dim}H^0(\Mod{d}{Y}; \Line{\rho \circ \beta } ) +
( \lambda (\rho ) + 1)^{g-1} \sum_{[Z] \neq 0} 
q (\sigma , \rho (\beta P), Z ) \right) 
\end{equation}
From what was said in the remark proceeding this proof we know
that 
\begin{align*}
\sum_{ [Z] } 
q (\sigma , \rho (\beta P) , Z ) & =
(-1)^{w_2 (\rho) \cdot j } \sum_{ [Z] }  
(-1)^{\epsilon (\sigma + Z) - \epsilon (\sigma )} \\
& =
(-1)^{w_2 (\rho \circ \beta P) + \epsilon (\sigma )} 
\sum_{ \sigma'  } (-1)^{\epsilon (\sigma')} \\
& =
(-1)^{w_2 (\rho \circ \beta P)   + \epsilon (\sigma )}
( (2^{2g-1} + 2^{g-1} ) - (2^{2g-1} - 2^{g-1} ) )  \\
& =
 (-1)^{w_2 (\rho \cdot \beta P)   + \epsilon (\sigma )} 2^g
\end{align*}
so that 
$$
\sum_{ [Z] \neq 0 } q (\sigma , \rho (\beta P), Z ) =
(-1)^{w_2 (\rho \circ \beta P) } ( (-1)^{\epsilon (\sigma )} 2^g - 1 ) .
 $$
Plugging this into \eqref{dim1} gives us the proposition. 
\end{proof}

\subsection{Application to spin Chern-Simons theory}

In this section we use Proposition \ref{dim form} to compute
the dimensions of the Hilbert spaces of the spin Chern-Simons
theory for closed, spin 2-manifolds of genera greater than one.  The 
genus one case is considered in the following chapter.  For the most
part, our job is already done.  However, some care must be taken 
to obtain the ``correct''  quantum theory.  In particular this requires
the proper choice of classical level.  Indeed, to obtain the quantum theory
which we conjecture to correspond with the spin-TQFT of Blanchet
and Masbaum, we must consider only a certain subset of all possible levels.
(Actually, as future work will show, there is a more comprehensive
$SO_3$ theory that incorporates all possible levels, and the spin 
Chern-Simons we consider here is a subset of this theory.)

To discuss the levels consider the standard representation 3-dimensional
representation $\text{id}_{SO_3}$ and define
$$
{\bf 1} = \lambda ( \text{id}_{SO_3} - 3) \in E^4(BSO_3)
$$
then ${\bf 1}$ generates the levels of $SO_3$ so that $E^4(BSO_3) = 
\Z \cdot {\bf 1}$.
We let the boldfaced ${\bf k}$ denote the level $k \cdot {\bf 1}$
for any $k \in \Z$.  To obtain the spin-Chern-Simons that we want
we consider only even valued levels.  According to the prescription
described in subsection \ref{2.2.1} and the geometric equivalence 
determined by proposition \ref{rtK}, the Hilbert spaces we must
consider is 
\begin{equation}\label{sCS Hilb}
\Hilb (Y, \sigma, {\bf k}, w_2) = 
H^0(\Mod{w_2}{Y}; \Line{\bf k + 1}{(\sigma)}).
\end{equation}

Having established the identity of our Hilbert space, the computation
of its dimension is a trivial corollary of proposition \ref{dim form}.
Recall from Section 2.2, that the Pfaffian line bundles $\Line{\rho}
{(\sigma)}$ are graded and that the grading is determined by the 
mod 2 index of the twisted Dirac operators.  As this is a topopological
invariant the grading is obviously constant over connected components
of the moduli stack.  In particular the grading of the line bundles over 
$\Mod{w_2}{Y}$ will depend on $w_2$ and the level. 
At level ${\bf k + 1}$, $k$ even, the line bundle 
has odd grading when $w_2 =1$ and has even grading when $w_2 = 0$. 

\begin{corollary}\label{dim Hilb}
Let $Y$ be a closed, oriented genus-$g$ 2-manifold with spin stucture 
$\sigma$.  Then the dimension of the even Hilbert space of the quantum 
spin-Chern-Simons theory at level ${\bf k}$ is given by
\begin{align}\label{dim even Hilb}
\dim \Hilb (Y, \sigma, {\bf k}, w_2=0) &   = \\ \nonumber
\frac1{2^{2g}} 
( \dim H^0(\Mod{d=0}{Y}; \Line{\bf 2k + 2} ) + 
 ( (-1)^{\epsilon (\sigma )} & 2^g -1) (k + 2)^{g-1})
\end{align}
and the dimension of the odd Hilbert space is given by
\begin{align}\label{dim odd Hilb}
\dim \Hilb (Y, \sigma, {\bf k}, w_2=1) &  = \\ \nonumber
\frac1{2^{2g}} 
( \dim H^0(\Mod{d=1}{Y}; \Line{\bf 2k + 2} ) - 
 ( (-1)^{\epsilon (\sigma )} & 2^g -1) (k + 2)^{g-1}).
\end{align}
\end{corollary}
This formula agrees with the one given in Theorem 19.1 in \cite{BM}.  

We consider, in particular, the even Hilbert spaces. The next proposition
relates the Hilbert spaces of our $SO_3$ quantum theory to the Hilbert space
of the well-known $SU_2$ quantum theory.  To do so we must first say a few 
words about the $SU_2$ quantum theory.  We start with the levels for the
$SU_2$ theory.

To discuss the levels consider the realization of the standard $\C^2$ 
representation $\rho: SU_2 \rightarrow SO(\R^4)$ and define
$$
{\bf 1'} = \lambda ( \rho - 4 ) \in E^4(BSU_2)
$$
then ${\bf 1'}$ generates the levels of $SU_2$ so that $E^4(BSU_2) = 
\Z \cdot {\bf 1'}$.
We let the boldfaced ${\bf k'}$ denote the level $k \cdot {\bf 1'}$
for any $k \in \Z$.  We point out that the adjoint representation $\text{ad} :
SU_2 \rightarrow SO(\Lie{so}_3)$ represents the level $\lambda 
(\text{ad}) = {\bf 2'}$.  The homomorphism $\beta : SU_2 
\rightarrow SO_3$ induces a homomorphism 
\begin{align*}
\beta^* : E^4 (BSO_3) & \longrightarrow E^4 (BSU_2) \\ 
               k \cdot {\bf 1} & \longmapsto 2k \cdot {\bf 1'}
\end{align*}
According to the prescription described in section \ref{2.2.1} and the 
geometric equivalence determined by proposition \ref{rtK}, the Hilbert space 
we must consider is 
\begin{equation}\label{CS Hilb}
\Hilb' (Y, {\bf k'}) = 
H^0(\Mod{d=0}{Y}; \Line{\bf k' + 2'}{(Y)}).
\end{equation}
It is easy to show that, over $\Mod{d=0}{Y}$,  
$\Line{\rho \circ \beta}$ does not depend on the spin structure
$\sigma$ in the sense that for any two spin structures there exists 
a connection preserving isometry between the two corresponding 
line bundles which is unique up to a factor in $\T$.  This is why $\sigma$
does not appear in the denotation for the $SU_2$ Hilbert space. 

The definitions \eqref{sCS Hilb} and \eqref{CS Hilb} imply
that we have inclusions 
$$
i_{\sigma} : \Hilb (Y, \sigma, {\bf k}, w_2=0) 
\hookrightarrow \Hilb' (Y, {\bf 2k'})
$$ 
and that each these subspaces is the image of the (respective) projection 
$$
P_{\sigma} = \frac{1}{2^g}
\sum_{Z \in \Mod{\zmod2}{Y}} [Z]_{CS, \sigma}.
$$
We offer the following proposition.

\begin{proposition}\label{projs}
Take any non-trivial element $\ell \in H^1(Y;\zmod2)$.  Then
$$
P_{\sigma + \ell} \circ P_{\sigma} = 0.
$$ 
\end{proposition}

\begin{proof}
We first compute $[Z]_{CS, \sigma + \ell}$ in terms of 
$[Z]_{CS, \sigma}$.  If we compose the former with the 
inverse of the latter we get a covariantly constant automorphism
of $\Line{\bf 2k'}$ which projects to the identity on $\Mod{SU_2}{Y}$.
This is just multiplication by some constant in $\T$ so that to compute
it we only need to do so over a fixed point of $[Z]$.  Thus
\begin{align*}
[Z]_{CS, \sigma + \ell} \circ [Z]_{CS, \sigma}^{-1} & =
q (\sigma + \ell, 0, Z) \cdot q (\sigma, 0, Z)^{-1} \\ 
& = (-1)^{\text{ind}_2(D_{\sigma + \ell + Z}) 
- \text{ind}_2(D_{\sigma + \ell})  
- \text{ind}_2(D_{\sigma + Z}) 
+ \text{ind}_2(D_{\sigma})} \\
& = (-1)^{Z \smile \ell}.
\end{align*}
where the last equality follows from the fact that 
$$
\text{ind}_2(D_{\sigma + Z}) - \text{ind}_2(D_{\sigma})
$$
is quadratic with respect to $Z$ and the corresponding bilinear form
on $H^1(Y;\zmod2)$ is the one determined by the cup product.
Finally we compute
\begin{align*}
P_{\sigma + \ell} \circ P_{\sigma} = &
\sum_{Z'} \sum_Z [Z']_{CS, \sigma + \ell}
 \circ [Z]_{CS, \sigma} \\
= & \sum_{Z'} \sum_Z (-1)^{Z' \smile \ell}[Z']_{CS, \sigma} 
\circ [Z]_{CS, \sigma} \\
= &  \sum_{Z'} \sum_Z (-1)^{Z' \smile \ell}
[Z' + Z]_{CS, \sigma} \\
= & \left( \sum_{Z'} (-1)^{Z' \smile \ell} \right) 
P_{\sigma} \\
= & 0
\end{align*}
\end{proof}

The upshot to this proposition is that the subspaces 
 $\Hilb (Y, \sigma, {\bf k}, w_2=0) \subset \Hilb' (Y, {\bf 2k'})$
 are disjoint for different spin structures.  In fact, that the projections
 $P_{\sigma}$ are constructed out of isometries of the Hilbert space
 implies that the subspaces are orthogonal to each other.  A straightforward
 dimension count shows that
 $$
 \sum_{\sigma} \dim \Hilb (Y, \sigma, {\bf k}, w_2=0) 
 = \dim \Hilb' (Y, {\bf 2k'}).
  $$
  
  We tie all of this together to conclude with the final proposition of this
  paper which is the Hamiltonian version of the Proposition \ref{formal refine}.
  
  \begin{proposition}\label{decomp}
  For $k$ even, we have the following orthogonal decompositon
  $$
  \Hilb' (Y, {\bf 2k'}) = \bigoplus_{\sigma}
  \Hilb (Y, \sigma, {\bf k}, w_2=0) 
  $$ 
  so that the Hilbert spaces for the $SO_3$ spin-Chern-Simons theory
  are refinements of the Hilbert space for the $SU_2$ theory.
  \end{proposition}

\bibliographystyle{plain}

\end{document}